\documentclass[a4paper,10pt,oneside,onecolumn,preprint,doubleblind]{article}
\usepackage{geometry,amssymb,hyperref}
\usepackage[many]{tcolorbox}

\newtheorem{te}{Theorem}
\newtheorem{lm}{Lemma}

\newenvironment{dkz}{{\bf Proof:}}{$\square$}
\newcommand{\torka}[1]{\langle{#1}\rangle}

\newcommand{\upar}[2]{\torka{{#1},{#2}}}
\newcommand{\izbaci}[1]{}

\newcommand{\Part}[1]{{\cal P}{(#1)}}
\newcommand{\dom}[1]{\mathrm{dom}({#1})}

\newcommand{\interp}[2]{\mathrm{int}_{#1}({#2})}

\newcommand{\funure}[3]{{\mathrm{Fn}}({#1},{#2},{#3})}

\newcommand{\ZFC}{\mathrm{ZFC}}
\newcommand{\ZF}{\mathrm{ZF}}

\newcommand{\jez}[1]{{\mathrm L}_{#1}}

\newcommand{\BPI}{{\mathrm{BPI}}}
\newcommand{\DDF}{{\mathrm{DDF}}}
\newcommand{\PG}{{\mathrm{PG}}}

\newcommand{\skupsvih}[2]{{\{{#1}\,|\,{#2}\}}}
\newcommand{\jedef}{:=}
\newcommand{\ekdef}{\stackrel{\mathrm{def}}\Leftrightarrow}

\newcommand{\HLT}{Halpern-L\" auchli theorem}
\author{Nedeljko Stefanovi\' c}
\title{The ZFC analogue of the Halpern-Levy theorem}
\begin{document}
	\maketitle		
		\begin{abstract}
			Here we present $\ZFC$ theorems yielding the \HLT{} and avoiding metamathematical notions in their formulations.
		\end{abstract}

	\section{Introduction}
	
	In order to prove the consistency of the theory $\ZF$ with the negation of the Axiom of Choice, Paul Cohen in \cite{cohen1} and \cite{cohen2} has constructed a model%
	\footnote{In the notation of this article, that model is $N(\omega,\funure{\omega}{\{0,1\}}{\aleph_0},M)$.}
	that has been a subject of intensive study since then, a model which is today known as the Cohen symmetric model.
	Not long afterwards, Halpern and Levy proved in \cite{hl2} that in the Cohen symmetric model the statement that every Boolean algebra has an ultrafilter is true.
	Their proof was crucially based on the \HLT{}, a Ramsey-theoretic statement about product of trees.
	This, today famous and deep, theorem  was proven in the paper \cite{hl} of Halpern and L\" auchli.
	Later, Harrington found another proof of the \HLT{} that uses forcing.
	Expositions of versions of Harrington's proof can be found in \cite{tf}, \cite{omar}, \cite{jtm}  and  \cite{dobrinen}.
	
	Since the original papers of Halpern-L\" auchli and Halpern-Levy, many other applications of the \HLT{} have appeared in the literature (see, for example,\cite{rs}).
	In the paper \cite{ns} the \HLT{} is derived from the fact that $\BPI$ holds in the Cohen symmetric model,
	and it is shown how that model can be used instead of the \HLT{} to prove statements which are otherwise proved
	by applying the theorem itself.
	More precisely, a new mathematical method was given in \cite{ns}, a method that can be seen as an alternative to the \HLT{} in its applications and a method that has the power not weaker than the \HLT{}.
	Moreover, our paper \cite{ns} gives a proof of relative consistency, with $\ZFC$, of a coloring principle $\DDF$ (introduced in \cite{z})  which is stronger than the \HLT{}. In fact, the paper \cite{ns} gives a proof of the relative consistency of the principle  $\DDF$  for a set of colors whose cardinality is less than $2^{\aleph_0}$, and with the continuum being a limit cardinal.	
	This provides a yet another  proof of the \HLT{} from the axioms of $\ZFC$.
	
	Not long afterwards,  Lambie-Hanson and Zucker in \cite{cz} and independently this author have proved a consistency result
	for another coloring principle $\PG$  (stronger than $\DDF$) about products of Polish spaces by at most countably many colours.
	This principle  $\PG$  easily implies the \HLT{}. Later, the author of this paper generalized \cite{nsa} this consistency result for a set of colors of cardinality less than $2^{\aleph_0}$
	in a model where the continuum is a limit cardinal.
	The proof presented in \cite{nsa} differs from the proof presented in \cite{cz}.
	In this paper, we also give some new $\ZFC$ theorem yielding the \HLT{} while avoiding metamathematical notions in their formulations.
	
	\subsection{Cohen's symmetric model and its generalizations}
	\vspace{1em}
	
	We shall need the following important fact saying that $M$ is a class in the model $M[G]$, definable with a parameter in $M$.
	
	\begin{te}{\upshape\cite[Theorem 2.16]{audrito}}
		{\upshape({\itshape Laver's Theorem})}
		\label{lejverovaTeorema}
		If a transitive model $M$ satisfies a suitable finite fragment of the theory $\ZFC$
		and if $P\in M$ is some atomless partial ordering, and $G$ is some $P$-generic filter over $M$,
		then $M$ is a class in the model $M[G]$, definable with a parameter in $M$.
	\end{te}
		
	We shall use the same definitions and notation as in \cite{ns}, but for the convenience of the reader, we reproduce some of them here.
	Let $M$ be a countable transitive model for a suitable finite fragment of the theory
	$\ZFC$, and let $I\in M$ be some infinite set. Let $A$, $B$ and $\kappa$ be such that
	\begin{equation}
		\label{ABKappa}
		\aleph_0\leqslant\kappa\leqslant|A|,\quad |B|\geqslant 2.
	\end{equation}
	Let us define the ordering $\funure AB\kappa$ as the set of all partial functions
	from $A$ to $B$ of cardinality less than $\kappa$ and ordered by reverse inclusion.
	If $\kappa=\aleph_0$, this is the set of all finite partial functions from $A$ to $B$.
	
	Let $Q$ denote some atomless partial order belonging to the model
	$M$, and let $P$ denote the partial ordering
	of all finite partial functions from $I$ to $Q$ ordered as follows:
	$$
	(\forall s_1,s_2\in P)(s_1\leqslant{s_2}\ekdef(\dom{s_2}\subseteq\dom{s_1}\,\land\,
	(\forall i\in\dom{s_2})s_1(i)\leqslant s_2(i))).
	$$
	
	Let $G$ denote some $P$-generic filter over $M$, and let $X$ denote
	some class in the model $M[G]$, definable by parameters in $M$ and such that $X\subseteq M$.
	Also, for each $i\in I$ we define
	
	$$
	\dot c_i\jedef\skupsvih{\upar p{(p(i))\check{}\,}}{
		p\in P,\,i\in\dom p},
	\quad c_i\jedef\interp G{\dot c_i}
	$$
	and
	$$
	\dot C\jedef\skupsvih{\upar p{\dot c_i}}{p\in P,\,i\in I},\quad C=\interp G{\dot C}.
	$$
	
	Obviously,
	\begin{equation}
		\label{cPrvo}
		c_i=\skupsvih{p(i)}{p\in G\land i\in\dom p},\quad c_i\subseteq Q\in M,
	\end{equation}
	holds for all $i\in I$, and also
	$$
	C=\skupsvih{c_i}{i\in I}.
	$$
	
	By the forcing product lemma
	\begin{equation}
		\label{nezavisnoGenericki}
		c_i\mbox{ is }Q\mbox{-generic filter over }M[f_i],\quad\mbox{where }f_i:I\setminus\{i\}\longrightarrow C,\,f_i(j)=c_j
	\end{equation}
	holds for every $i\in I$.
	In particular, for any distinct $i_1,i_2\in I$ there exist $s_1\in c_{i_1}$ and $s_2\in c_{i_2}$
	such that $s_1\bot s_2$,
	and therefore $c_i\neq c_j$, so the set $C$ is infinite.
	Let us define the model	
	$$
	N(I,Q,X)\jedef HOD(C\cup\{C\}\cup X)^{M[G]}.
	$$
	
	In this definition, we ignore the small inaccuracy that $G$ is not uniquely determined by $I$, $Q$, $M$ and $X$.
	If there is no danger of confusion, instead of $N(I,Q,X)$ we will write $N$.
	
	Cohen's symmetric model typically means the model
	$N(\omega,\funure\omega{\{0,1\}}{\aleph_0},X)$, where $X=\emptyset$ or $X=M$.
	We will work under the assumption that $Q\subseteq N$. Then
	$(\forall i\in I)c_i\in N$, and therefore $C\in N$. The assumption is
	satisfied when $X=M$ because then
	$Q\subseteq M\subseteq N$.
	
	\subsection{Properties of  Cohen's symmetric model}
	\vspace{1em}
	
	Here we use a more general approach of introducing Cohen's symmetric model.
	It is not only about adding countably many Cohen reals, but about adding arbitrarily many generic objects of any type.
	A neighborhood is by definition any set of the form
	$$
	[s]\jedef\skupsvih{c\in C}{s\in c},\quad s\in Q.
	$$
	
	A neighborhood of an element $c\in C$ is by definition any neighborhood to which $c$ belongs.
	For any $s\in Q$ and $r\in C$:
	\begin{equation}
		\label{karakterizacijaOkoline}
		r\in[s]\Leftrightarrow s\in r.
	\end{equation}
	
	For any incompatible $s_1,s_2\in Q$ we have $[s_1]\cap[s_2]=\emptyset$.
	Therefore different elements of the set $C$ can be placed in disjoint neighborhoods.
	
	Let us now turn to the properties of the Cohen's symmetric model. The following theorem highlights the importance
	of one simple type of automorphism of the ordering $P$.
	\begin{te}{\upshape\cite[Theorem 2]{ns}}
		\label{lepiAutomorfizam}
		Let $\sigma\in M$ be a permutation of the set $I$, and let $f$ be
		an automorphism of the ordering $P$ such that $f(x)=x\circ\sigma$ for any
		$x\in P$.
		Then $f\in M$, and $H=f[G]$ is a $P$-generic filter over $M$, and $(\forall i\in I)\interp H{\dot c_i}=c_{\sigma(i)}$, $M[H]=M[G]$ and $\dot C_H=C$.
	\end{te}
	\vspace{1em}
	
	The most important property of the Cohen's symmetric model is given by the continuity lemma.
	Before that, we give a similar theorem for the model $M[G]$.	
	Note that in the following theorem, the parameters used to define the class $X$ need not be
	a subset of $\{a_1, ...,a_m\}$.
	
	\begin{te}{\upshape({\itshape Continuity}) \cite{omar}, \cite[Theorem 4]{ns}}
		\label{lemaONeprekidnosti}
		Assume the following holds:
		\begin{enumerate}
			\upshape
			\item {\itshape  $a_1,\dots,a_m\in M$,}
			\item {\itshape $r_1,\dots,r_n$ are mutually distinct elements of the set $C$,}
			\item {\itshape $\varphi(x_1,\dots,x_m,y_1,\dots,y_n,z)$ is the formula of the language $\jez{\ZFC}$,}
			\item {\itshape In the model $M$, some theorem of the theory $\ZFC$ (depending only on $\varphi$) is true.}
		\end{enumerate}
		
		Then there are disjoint neighborhoods
		$[s_1],\dots,[s_n]$ of elements $r_1,\dots,r_n$ respectively, such that for all
		$r_1'\in[s_1],\dots,r_n'\in[s_n]$
		$$
		N\models(\varphi(a_1,\dots,a_m,r_1,\dots,r_n,C)\Leftrightarrow
		\varphi(a_1,\dots,a_m,r_1',\dots,r_n',C)).
		$$
		$$
		M[G]\models(\varphi(a_1,\dots,a_m,r_1,\dots,r_n,C)\Leftrightarrow
		\varphi(a_1,\dots,a_m,r_1',\dots,r_n',C)).
		$$
	\end{te}
	\begin{te}{\upshape\cite[Theorem 6]{ns}}
		\label{najmanjiSkupParametara}
		If $X\subseteq M$, then for each $b\in OD(X\cup C\cup\{C\})^{M[G]}$
		there is the inclusion-smallest set
		$A\subseteq C$ such that $b\in OD(X\cup A\cup\{C\})^{M[G]}$.
		In addition, the set $A$ is finite.
	\end{te}
	
	We will further denote by $F_1$ the function 
	which maps each $b\in N$ to
	the inclusion-smallest subset $A$ of the set $C$ such that $b\in OD(X\cup A\cup\{C\})^{M[G]}$.
	
	\begin{te}{\upshape\cite[Theorem 7]{ns}}
		\label{F1Klasa}
		The function $F_1$ is a class in the model
		$M[G]$ definable with parameters from the set
		$$
		ORD^M\cup\{C\}\cup T,
		$$
		where $T$ is the set of parameters by which the class $X$ is definable in $M[G]$.
	\end{te}
	\vspace{1em}
	
	The Cohen's symmetric model satisfies the axioms of $\ZF$, but not the Axiom of Choice.
	
	In \cite{ns} is given the theorem of Halpern and Levy \cite{hl2} for the generalized version of the Cohen's symmetric model.
	
	\begin{te} {\upshape\cite[Theorem 16]{ns}}
		\label{bpi}
		Assume that $X$ is a class of $M[G]$ definable by parameters
		from the set
		$$
		X\cup\{C\}\cup ORD^M
		$$
		and that for every $\alpha\in ORD^M$ there is a well ordering of the
		set $X\cap V_\alpha^N$, which belongs to the model $N$.
		Then in the model $N$ every Boolean algebra has a prime ideal.
	\end{te}
		
	Polish space is a complete separable metric space. Perfect space is a Polish space without isolated points.
	
	\begin{te}{\upshape\cite[Theorem 5]{nsa}}
		\label{gusto}
		Every perfect space contains a dense $G_\delta$ subspace homeomorphic to the space ${}^\omega\omega$.
	\end{te}
		
	\section{The results} 
	
	We start with the following simple fact.
	
	\begin{te}\label{BerZaBerove}
		Let ${\cal P}_1=\upar{P_1}{\rho_1},\dots,{\cal P}_d=\upar{P_d}{\rho_d}$ be Polish spaces,
		$M_1,\dots,M_d$ subsets of $P_1,\dots,P_d,$ respectively, all of the being of the first category and 
		$f:P_1\times\cdots\times P_d\longrightarrow\omega$ such that for each $n\in\omega$
		the set $f^{-1}[\{n\}]$ has the Baire property. Then there are some dense subsets $D_1,\dots,D_d$
		of spaces ${\cal P}_1,\dots,{\cal P}_d,$ respectively, such that $\bigwedge_{i=1}^d D_i\cap M_i=\emptyset$ and the function $f$ is locally constant on the set $D_1\times\cdots\times D_d$ .
	\end{te}
	\begin{dkz}
		Note that isolated point (if it exists) is somewhere dense set. Therefore, without loss of generality we can assume that spaces
		${\cal P}_1=\upar{P_1}{\rho_1},\dots,{\cal P}_d=\upar{P_d}{\rho_d}$ have no isolated points.
		Also, without loss of generality we can assume that for each $i$ the set $M_i$ is countable union of closed nowhere dense sets.
		The  set $P_i\setminus M_i$ forms a Polish space in relative topology with some metric for each $i$. Therefore, without loss of generality we can assume that
		$M_i=\emptyset$ for each $i$.
		According to Theorem \ref{gusto}, without loss of generality we can assume that ${\cal P}_i$ is the space ${}^\omega\omega$ for each $i$.
		
		Let $U_n$ be an open set in the space and $M_n'$ a countable union of nowhere dense closed sets in the space ${\cal P}_1\times\cdots\times{\cal P}_d$
		with $U_n\triangle f^{-1}[\{n\}]\subseteq M_n'$ for each $n$.  Let us choose the countable transitive model $M$ containing the codes for sequences $U=(U_n)_{n\in\omega}$ and $M'=(M_n')_{n\in\omega}$ for all $n$. 
		
		We will use the notation given in the introduction. Let $I=\omega$ and $Q=\funure\omega\omega{\aleph_0}$
		and $x_i=\bigcup c_i$ for each $i\in I$. Let  $D=D_1\times\cdots\times D_d$ where $D_i=\skupsvih{x_i}{i\in I}$ for each $i\in\{1,\dots,d\}$.
		
		Let $\torka{t_1,\dots,t_d}\in D$ be arbitrary. Cohen real can't be in a set of the first category coded in $M$. Therefore,
		for $n=f(t_1,\dots,t_d)$ it holds $\torka{t_1,\dots,t_d}\in U_n\cap M[G]$. By the assumption that the set $U_n$ is open
		there are finite sets $s_1\subseteq t_1,\dots,,s_d\subseteq t_d$ such that  $[s_1]\times\cdots\times[s_d]\subseteq U_n$ .
		All members of the set $[s_i]$ are Cohen reals and therefore $f$ is constantly equal to $n$ on the set $[s_1]\times\cdots\times[s_d]$.
	\end{dkz}
	
	Unfortunately, it is not possible to derive the \HLT{} from this lemma. For this reason we give another theorem
	which easily implies the \HLT{}, but we need some facts before its statement. 
	
	We will deal with rooted trees where every node is at a finite height with at least one but at most
	countably many immediate successors with usual tree topology on the set of all branches. 
	
	We say that the branch is condensational if every neighborhood of $b$ is uncountable.
	We say that the node is condensational if it is member of at least one condensational branch.
	The set of all condensational nodes is either empty or it forms a subtree (which means that if at least one of the
	immediate successors of  some node $v$ is condesational  then the node $v$ is condensational).
	We define the condensational subtree as subtree of all condensational nodes.
	
	Of course, a branch is condensational iff each of its nodes is condensational.
	The set of all condensational branches is equal to set of all branches of the condensational subtree
	and the set of condensational nodes is equal to the set of all nodes of the condensational subtree.
	A tree is equal to it's condensational subtree iff it has no isolated branches.
	
	\begin{lm}
		These notions are absolute for transitive models of sufficiently large fragment of the $\ZF$ theory.
	\end{lm}
	\begin{dkz}
		For that purpose we will slightly generalize the notion of tree
		by removing the condition that every node has at least one immediate successor.
		
		For the tree $T$ we define it's derivative $T'$ as the subtree of the tree $T$ where the set of nodes of the tree $T'$
		is the set of nodes of the tree $T$ belonging to more than one branch of the tree $T$.
		It is easy to check that $T'$ is subtree of the tree $T$.
		
		The tree $T$ has no isolated branches iff $T'=T$. The cardinal number of the set $[T]\setminus[T']$ of branches of the tree $T$ which are not branches of the tree $T'$
		is not greater than the cardinal number of the set of nodes of the tree $T$ which are not nodes of the tree $T'$. 
		The set of nodes is at most countable. Therefore the set $[T]\setminus[T']$ is at most countable.
		
		For a given tree $T$ we can construct the sequence $(T_\alpha)_{\alpha<\omega_1}$ of trees where $T_0=T$,
		$T_{\alpha+1}=T_\alpha'$ and for infinite limit ordinal $\alpha<\omega$ the tree $T_\alpha$ is the subtree
		of the tree $T$ with the set of nodes equal to the intersection of the set of nodes of trees $T_\beta$ for all $\beta<\alpha$.
		
		For each $\alpha<\omega$ the set $[T]\setminus[T_\alpha]$ is at most countable. Therefore
		every neighborhood of any condensational branch $b$ of the tree $T$ has uncountably many branches
		in the subtree $T_\alpha$. Therefore every condensational node of the tree $T$ is also a node of the tree $T_\alpha$.
		
		The set of nodes is at most countable. Therefore $T_{\alpha+1}=T_\alpha$ for some $\alpha<\omega_1$. 
		If $T_\alpha$ is an empty tree then the tree $T$ has at most countably many branches, has no condensational
		branch and has no condensational node. If $T_\alpha$ is nonempty, then it is a tree without isolated branches.
		Therefore every node of the tree $T_\alpha$ is condensational and $T_\alpha$ is condensational subtree of
		the tree $T$.
		This construction is absolute for all transitive models of some theorem of $\ZF$.
	\end{dkz}
	\vspace{1em}
	
	Hereafter we will limit our considerations to trees without isolated branches. Also, we will
	consider trees with some linear ordering of the set of nodes where the set of immediate successors
	of any node is either finite or has order type $\omega$. It gives us one linear ordering of branches.
	
	We say that branches $b_1$ and $b_2$ are near if $b_1=b_2$ or $b_1\neq b_2$
	with no branches between $b_1$ and $b_2$. This relation is obviously reflexive and symmetric.
	It is also transitive because of the nonexistence of isolated branches.
	
	After the identification of near branches we obtain a dense linear ordering. After removing
	the minimum and the maximum (if they exists) we obtain a dense linear ordering without endpoints. 
	Let us denote it's domain by $B$.
	
	To each node except the root we can associate the leftmost branch containing that node. The set $D$ of all such branches (assuming the identification of near branches) is a countable dense set.
	
	By the Cantor's back and forth procedure we can define an increasing bijection $f$ between the set $\mathbb Q\cap(0,1)$ and the set $D$. 
	It is possible to define a construction of such mapping $f$ so that the construction is absolute for transitive models of some fragment $\ZF$
	containing that tree together with the ordering of nodes.
	
	Let us prove that $f$ can be extended to exactly one increasing mapping $g:(0,1)\longrightarrow B$.
	For any $a\in(0,1)\setminus\mathbb Q$ let
	$$
	L(a)=\skupsvih{f(x)}{x\in(0,1)\cap\mathbb Q,\,x<a},\quad R(a)=\skupsvih{f(x)}{x\in(0,1)\cap\mathbb Q,\,x>a}.
	$$
	
	Then $L(a)<R(a)$ and $L(a)\cup R(a)$ is equal to the dense set $D$. Therefore there are no different branches
	$b_1$ and $b_2$ between sets $L(a)$ and $R(a)$. It is easy to define the leftmost branch $b$ with $L(a)<b$.
	Then we define $g(a)$ as the branch $b$. 	It is easy to check that $g$ is an increasing extension of the function $f$.
	For any $b\in B\setminus D$ we can find the unique irrational number $a$ with
	$$
	\skupsvih{x\in(0,1)\cap\mathbb Q}{f(x)<b}<a<\skupsvih{x\in(0,1)\cap\mathbb Q}{f(x)>b}.
	$$
	
	Then $g(a)=b$ and therefore $g$ is a surjection onto the set $B$. All steps are absolute for transitive models
	of some $\ZF$ theorem.
	
	If we identify a subset of $\omega$ with an infinite sequence of zeros and ones, then we can define topology on $\Part\omega$ as the 
	product topology of ${}^\omega2$. Then we can define the notion of $F_\sigma$ filter as a filter over $\omega$
	which is a countable union of closed sets.
	
	Note that in the used topology of $\Part\omega$, the union and the intersection of a pair of subsets of $\omega$ are continuous operations
	and that the space $\Part\omega$ is a compact metrizable space. For every $k$ the mapping $f_k$ defined by
	$$
	f_k:\Part\omega^{k+1}\longrightarrow\Part\omega,\quad f_k(A_1,\dots,A_k,B)=(A_1\cap\cdots\cap A_d)\cup B
	$$
	is continuous. Let $F$ be a filter generated by the union of the countable set $\skupsvih{X_k}{k\in\omega}$ of closed sets.
	Sets $X_k$ are compact as closed sets in a compact Hausdorff space. 
	Finally, $F$ is an $F_\sigma$ filter as a countable union of compacts sets  $f_k[X_1\times\cdots\times X_k\times\Part\omega]$
	for all $k\in\omega$.
	
	In particular countably generated filters are $F_\sigma$ filters. Also, for every at most countable family 
	of $F_\sigma$ filters the filter generated by their union (if it exists) is also an $F_\sigma$ filter.
	Of course for every at most countable family of countably generated filters the filter generated by
	their union (if it exists) is a countably generated filter.
	
	Let us give an example of an $F_\sigma$ filter which is not countably generated. Let
	$$
	F'=\skupsvih{S\subseteq\omega}{\sum_{n\in\omega\setminus S}\frac 1{n+1}<\infty}.
	$$
	
	For $n\in\omega$ the set $Y(n)=\skupsvih{S\subseteq\omega}{n\in S}$ is clopen.
	For $m,k\in\omega$ the set $Z(m,k)=\skupsvih{S\subseteq\omega}{\sum_{n\in(\omega\setminus S)\cap k}\frac 1{n+1}\leqslant m}$
	is clopen as a Boolean combination of clopen sets $Y(n)$ for $n<k$.
	The filter $F'$ is an $F_\sigma$ filter because $F'=\bigcup_{m\in\omega}\bigcap_{k\in\omega}Y(m,k)$.
	
	Let us prove that the filter $F'$ is not countably generated. We will derive the contradiction from the assumption that the filter $F'$
	is generated by the set $\skupsvih{B_n}{n\in\omega}$ where $(\forall n)B_{n+1}\subseteq B_n$.
	
	The filter $F'$ has no finite elements and therefore $B_n$ is an infinite set for every $n$.	
	Let us choose a sequence $(e_k)_{k\in\omega}$ where $e_k\in B_k$ and $e_k\geqslant 2^k$ for all $k$.
	The set $S=\omega\setminus\skupsvih{e_k}{k\in\omega}$ belongs to the filter $F'$ but it is not superset of $B_k$ for any $k$.
	
	\begin{te}\label{uopstenjeHLT}
		Let $T_1,\dots,T_d$ be trees of height $\omega$ where every node has at least one but at most countably many immediate successors
		and $f:\bigcup_{l\in\omega}(T_1(l)\times\cdots\times T_d(l))\longrightarrow A$ where $A$ is a finite set with the discrete topology.
		Let $F_0$ be an $F_\sigma$ filter over $\omega$ and $M_1,\dots,M_d$ be subsets of the first category of spaces $[T_1],\dots,[T_d]$.
		
		Then there are dense sets $D_1,\dots,D_d$ in spaces $[T_1],\dots,[T_d]$ and $F_\sigma$ filter $F\supseteq F_0$ over $\omega$
		such that for each $i$, $D_i\cap M_i=\emptyset$ and that the mapping $g:D_1\times\cdots\times D_d\longrightarrow A$ defined as
		$$
		g(b_1,\dots,b_d)=a\quad\Leftrightarrow\quad\skupsvih{n\in\omega}{f(b_1(n),\dots,b_d(n))=a}\in F
		$$
		is continuous. If $F_0$ is a countably generated filter then the filter $F$ can be chosen as countably generated.
		Moreover $F$ can be chosen as a filter generated by $F_0\cup S$ for some countable set $S$.
	\end{te}
	\begin{dkz}
		Let $M$ be a countable transitive model containing $f,T_1,\dots,T_d$, codes for $M_1,\dots,M_d$, linear orderings of sets of nodes of trees $T_1,\dots,T_d$ so that the set of immediate successor of any node is finite or of order type $\omega$
		and so that the model $M$ contains a countable basis of $F_0$ if it is countably generated or Borel code of $F_0$ otherwise.
		
		We will use notation from the introduction. Let $N=N(\omega,\funure\omega 2{\aleph_0},M)$ be Cohen symmetric
		model with added countably many Cohen reals over $M$ and with parameters in $M$. We will identify
		a branch with the corresponding real number.  Let $F_1=F_0\cap N$. The set $F_1$ is a filter over $\omega$ in $N$.
		
		Since $N\models\BPI$ we can choose an ultrafilter $F_2$ over $\omega$ in $N$ with $F_2\supseteq F_1$.
		Let $C_0$ be a finite subset of $C$ such that $F_2$ is definable in $N$ by parameters from $C_0\cup M$.
		Let $I_0$ be a finite subset of $\omega$ with $(\forall i\in\omega)(i\in I_0\Leftrightarrow r_i\in C_0)$.
		Let us choose disjoint infinite subsets $I_1,\dots,I_d$ of $\omega\setminus I_0$ in $M$.
		
		For each $k\in\{1,\dots,d\}$ let $D_k$ be the set of branches of the tree $T_k$ corresponding to Cohen reals from the set $I_k$. For each $k$ the set $D_k$ is dense and disjoint with $M_k$. Let
		$$
		J=\skupsvih{\torka{c_1,\dots,c_d}\in(C\setminus C_0)^d}{\bigwedge_{i<j}c_i\neq c)j}.
		$$
		and let $D$ be the set of all $\torka{b_1,\dots,b_d}$ where $b_1,\dots,b_d$ are branches of trees
		$T_1,\dots,T_d$ corresponding to Cohen reals $c_1,\dots,c_d$ with $\torka{c_1,\dots,c_d}\in J$.
		Let us note that the set $D$ belongs to the model $N$.
		By Theorem \ref{lemaONeprekidnosti} the mapping
		$$
		h:D\longrightarrow A,\quad h(b_1,\dots,b_d)=a\,\Leftrightarrow\,\skupsvih{n\in\omega}{f(b_1(n),\dots,b_d(n))=a}\in F_2
		$$
		is continuous as well as its restriction $g$ to the set $D_1\times\cdots\times D_d$.
		
		Let us prove that the set $F_0\cup F_2$
		generates a filter. Otherwise there are sets $X\in F_0$ and $Y\in F_2$ with $X\cap Y=\emptyset$.
		Then $X\subseteq\omega\setminus Y$ and therefore $\omega\setminus Y$ belongs to
		filters $F_0$ (because $X\in F_0$), $F_1$ (because $Y\in N$) and $F_2$ (because $F_2\supseteq F_1$). This is in contradiction with $Y\in F_2$.
	\end{dkz}
	\vspace{1em}	
	
	The derivation of the \HLT{} from Theorem \ref{uopstenjeHLT} is quite similar to the derivation from the $PG$ principle.
	However this is not a proof of the \HLT{} because we are using the Halpern-Levy theorem which is derived from
	the \HLT{}.
	
	By replacing the Cohen forcing with some other forcing notion it is possible to obtain similar results. We give the formulation
	for the random forcing.
	
	\begin{te}\label{saRandomom}
		Let $T_1,\dots,T_d$ be trees of height $\omega$ such that every node has at least one but at most countably many immediate successors
		and $f:\bigcup_{l\in\omega}(T_1(l)\times\cdots\times T_d(l))\longrightarrow A$ where $A$ is a finite set with the discrete topology.
		Let $F_0$ be an $F_\sigma$ filter over $\omega$ and let $M_1,\dots,M_d$ be sets of measure zero in spaces $[T_1],\dots,[T_d]$
		with the usual measure on trees. 
		
		Then there are dense sets $D_1,\dots,D_d$ in spaces $[T_1],\dots,[T_d]$ and an $F_\sigma$ filter $F\supseteq F_0$ over $\omega$
		such that for each $i$, $D_i\cap M_i=\emptyset$ and that the following map is Borel $g:[T_1]\times\cdots\times[T_d]\longrightarrow A$ where
		$$
		(\forall b_1\in D_1,\dots,b_d\in D_d,a\in A)(g(b_1,\dots,b_d)=a\quad\Leftrightarrow\quad\skupsvih{n\in\omega}{f(b_1(n),\dots,b_d(n))=a}\in F).
		$$
		If $F_0$ is a countably generated filter then the filter $F$ can be chosen to be countably generated.
		Moreover $F$ can be chosen as a filter generated by $F_0\cup S$ for some countable set $S$.
	\end{te}
	
	The proof of Theorem \ref{saRandomom} is analogous to the proof of Theorem \ref{uopstenjeHLT}, but there is no easy derivation 
	of the \HLT{} from Theorem \ref{saRandomom}.
	
	\vspace{5mm}
	


\begin{thebibliography}{XXX}
		\bibitem{audrito}{\scshape Audrito Giorgio}, {\itshape Characterization of Set-Generic Extensions},
		{\upshape Universit\' a Degli studi di Torino, facolt\' a di scienze m.f.n.;
			Corso di studi in matematica, (2010).}
		\\ {\upshape\url{http://www.logicatorino.altervista.org/matteo_viale/audrito.pdf}},
		\bibitem{cohen1}{\scshape Cohen, Paul J.}
		{\itshape The independence of the continuum hypothesis I},
		{Proceedings of the U.S. National Academy of Sciemces, 50: 1143-48.}, (1963).
		\bibitem{cohen2}{\scshape Cohen, Paul J.}
		{\itshape The independence of the continuum hypothesis II},
		{Proceedings of the U.S. National Academy of Sciemces, 51: 105-48.}, (1964).
		\bibitem{omar}{\scshape De La Cruz Omar}, {\itshape Halpern and Levy Redux} --
		{An unpublished manuscript based on a series of meetings with S. Todor\v cevi\' c at CRM Bellatera center near Barcelona, (2003).}
		\bibitem{dobrinen}{\scshape Dobrinen, Natasha} {\itshape Forcing in Ramsey theory},
		{Proceedings of the 2016 RIMS Symposium on Infinite Combinatorics and Forcing Theory}, (2017).
		\bibitem{hl}{\scshape Halpern, James D.; L\" auchli Hans}, {\itshape A partition theorem},
		{Transactions of the American Mathematical Society, 124: 360–367}, (1966).
		\bibitem{hl2}{\scshape Halpern James D., L\' evy Azriel},
		{\itshape The  Boolean  Prime  Ideal  Theorem  does  not  imply  the Axiom of Choice},
		Axiomatic Set Theory, Proceedings of Symposia in Pure Mathematics, vol. XIII,Part I,
		pp. 83–134, AMS, Providence, (1971).
		\bibitem{cz}{\scshape Lambie-Hanson Chris and Zucker Andy}, {\itshape Polish space partition principles and the Halpern-L\" auchli theorem}, {\upshape Submitted} (2022).
		\bibitem{jtm}{\scshape Tatch Moore Justin}, {\itshape The Method of Forcing,} {\upshape \url{https://arxiv.org/abs/1902.03235}} (2019).
		\bibitem{ns}{\scshape Stefanovi\' c, N.}, {\itshape Alternatives to the Halpern-L\" auchli theorem}; {\upshape Annals of the Pure and Applied Logic, Vol. 174,  Iss. 9, 103313, (2023)},
		\bibitem{nsa}{\scshape Stefanovi\' c, N.}, {\itshape The coloring principle for the product of polish spaces and the Halpern and L\" auchli's theorem};
		{\url{https://arxiv.org/abs/2209.07768}, (2022)},
		\bibitem{rs}{\scshape Todor\v cevi\' c Stevo}, {\itshape Introduction to Ramsey Spaces},
		{\upshape Princeton University Press,} (2010).
		\bibitem{tf}{\scshape Todor\v cevi\' c Stevo, Farah Ilijas}, {\itshape Some applications of the method of forcing}, Yenisei, Moscow, (1995).
		\bibitem{z}{\scshape Andy Zucker}, {\itshape A new proof of the 2-dimensional Halpern-L\"auchli theorem;}, (2017), {\upshape \url{https://www.math.cmu.edu/~andrewz/HL2d.pdf}}.
	\end{thebibliography}
\end{document}